\documentclass[a4paper, 10pt]{amsart}

\usepackage{etoolbox}
\makeatletter
\let\ams@starttoc\@starttoc
\makeatother
\usepackage[parfill]{parskip}
\makeatletter
\let\@starttoc\ams@starttoc
\patchcmd{\@starttoc}{\makeatletter}{\makeatletter\parskip\z@}{}{}
\makeatother
\usepackage{setspace}

\usepackage{color,soul}
\definecolor{red}{rgb}{1,0,0}
\setulcolor{red}

\usepackage{amssymb}
\usepackage{amsmath}
\usepackage{amsthm}
\usepackage{fancyhdr}
\usepackage{esint}
\usepackage{enumerate}

\usepackage[colorlinks=true,linkcolor=blue,citecolor=blue,urlcolor=blue]{hyper ref}
\usepackage{bookmark}

\swapnumbers
\newtheorem{lemma}{Lemma}[section]
\newtheorem{prop}[lemma]{Proposition}
\newtheorem{thm}[lemma]{Theorem}
\newtheorem{cor}[lemma]{Corollary}

\theoremstyle{definition}
\newtheorem{defn}[lemma]{Definition}
\newtheorem{example}[lemma]{Example}
\newtheorem{rem}[lemma]{Remark}
\newtheorem{ass}[lemma]{Assumption}

\numberwithin{equation}{section}

\renewcommand{\(}{\left(}
\renewcommand{\)}{\right)}
\renewcommand{\~}{\tilde}
\renewcommand{\-}{\bar}

\newcommand{\cn}{\colon}
\newcommand{\8}{\infty}

\newcommand{\R}{\mathbb{R}}

\renewcommand{\S}{\mathbb{S}}
\renewcommand{\H}{\mathbb{H}}

\renewcommand{\a}{\alpha}

\newcommand{\g}{\gamma}
\renewcommand{\d}{\delta}
\newcommand{\e}{\epsilon}
\renewcommand{\k}{\kappa}
\renewcommand{\l}{\lambda}

\newcommand{\s}{\sigma}
\newcommand{\p}{\varphi}
\newcommand{\G}{\Gamma}

\newcommand{\T}{\Theta}

\newcommand{\del}{\partial}

\newcommand{\inpr}[2]{\langle #1,#2 \rangle}
\newcommand{\fr}[2]{\frac{#1}{#2}}

\DeclareMathOperator{\osc}{osc}

\DeclareMathOperator{\dist}{dist}

\newcommand{\Thm}{\begin{thm}}
\newcommand{\eThm}{\end{thm}}
\newcommand{\Def}{\begin{defn}}
\newcommand{\eDef}{\end{defn}}
\newcommand{\Prop}{\begin{prop}}
\newcommand{\eProp}{\end{prop}}
\newcommand{\Rem}{\begin{rem}}
\newcommand{\eRem}{\end{rem}}
\newcommand{\Lem}{\begin{lemma}}
\newcommand{\eLem}{\end{lemma}}
\newcommand{\eq}{\begin{equation}}
\newcommand{\eeq}{\end{equation}}
\newcommand{\Ex}{\begin{example}}
\newcommand{\eEx}{\end{example}}
\newcommand{\pf}{\begin{proof}}
\newcommand{\epf}{\end{proof}}
\newcommand{\Cor}{\begin{cor}}
\newcommand{\eCor}{\end{cor}}
\newcommand{\Ass}{\begin{ass}}
\newcommand{\eAss}{\end{ass}}
\newcommand{\SAl}{\begin{align}\begin{split}}

\newcommand{\ra}{\rightarrow}
\newcommand{\hra}{\hookrightarrow}

\newcommand{\mc}{\mathcal}

\newcommand{\mrm}{\mathrm}

\newcommand{\hp}{\hphantom}

\protected\def\ignorethis#1\endignorethis{}
\let\endignorethis\relax
\def\TOCstop{\addtocontents{toc}{\ignorethis}}
\def\TOCstart{\addtocontents{toc}{\endignorethis}}

\begin{document}
\title[Pinching and asymptotical roundness for inverse curvature flows]{Pinching and asymptotical roundness for inverse curvature flows in Euclidean space}
\author{Julian Scheuer}
\address{Dr.~Julian Scheuer, Albert-Ludwigs-Universit{\"a}t, Mathematisches Institut, Abteilung f\"ur Reine Mathematik, Eckerstr.~1, 79104 Freiburg, Germany}
\email{julian.scheuer@math.uni-freiburg.de}

\thanks{During the preparation of this work the author was being supported by the DFG}
\date{\today}

\begin{abstract}
We consider inverse curvature flows in the $(n+1)$-dimensional Euclidean space, $n\geq 2,$ expanding by arbitrary negative powers of a 1-homogeneous, monotone curvature function $F$ with some concavity properties. We obtain asymptotical roundness, meaning that circumradius minus inradius of the flow hypersurfaces decays to zero and that the flow becomes close to a flow of spheres. 

\end{abstract}

\maketitle

\tableofcontents

%%%%%%%%%%%%%%%%%%%%%%%%%%%%%%%%%%%%%%%%%%%%%%%%%%%%%%%%%%%%%%%%%%%%%%%%%
\section{Introduction}
We consider inverse curvature flows in Euclidean space $\R^{n+1},$ $n\geq 2,$
\eq \label{flow} \dot{x}=F^{-p}\nu,\ 0<p<\infty,\eeq
where $F$ is a symmetric, monotone, homogeneous of degree 1 and concave curvature function, which is defined in an open, convex cone of $\R^{n},$ such that $F$ vanishes on its boundary, where in case $p>1$ we assume $\G=\G_{+}.$ Here $\nu$ is the outward normal to the flow hypersurfaces  of the flow
\eq x\colon [0,T^{*})\times M\ra \R^{n+1}\eeq with starshaped initial hypersurface $M_{0},$ where in case $p>1$ we assume that $M_{0}$ is strictly convex.

In \cite{Gerhardt:01/2014} the same flow was considered. Here it is shown for any $0<p<\infty$ that the maximal time of existence $T^{*}$ is characterized by the property 
\eq \inf |x|\ra\infty,\ t\ra T^{*},\eeq
and that the rescaled surfaces 
\eq \~{M}_{t}=\T^{-1}M_{t},\eeq 
where $\T=\T(t)$ is the radius of a suitable expanding sphere, converge to the unit sphere in $C^{\infty},$ cf. \cite[Thm.~1.1, Thm.~1.2]{Gerhardt:01/2014}.

The goal of this paper is the improvement of the asymptotical behavior of the flow. Let us explain our motivation to do this. To our knowledge, almost all of the existing results on classical smooth inverse curvature flows, cf. the end of this introduction for an overview, assume that the initial hypersurface may be written as a graph over a sphere. Of course, there are many possibilities to do this. For example take a sphere $\mc{S}$ centered at $q\in\R^{n+1}$ as initial hypersurface of the flow. It follows, that is evolves through expanding spheres centered at $q$ as well. Now take another point $z\neq q$ in the interior of the ball enclosed by $\mc{S}$ and write $\mc{S}$ as a graph over a sphere $\mc{S}_{z}$ around $z,$
 \eq \mc{S}=\{(u(x),x)\cn x\in \mc{S}_{z}\}.\eeq 
 Let $\mc{S}$ evolve and let $u(t,\cdot)$ denote the corresponding graph functions over $\mc{S}_{z}.$
 From the previous observations it is clear that the oscillation of $u$ can not decay to zero as $t\ra T^{*},$ even though circumradius minus inradius of the hypersurfaces \textit{is} constantly zero. Thus the choice of the sphere $\mc{S}_{z}$ is not optimal and $u$ does not reflect the nice spherical shape of the evolving surfaces. The optimal sphere would be one around $q.$ In this paper we are going to show that such an optimal sphere exists in the sense that the flow hypersurfaces will fit to a flow of spheres arbitrarily close. We are going to achieve this for $0<p<\8$ and for $n\geq 2.$ The main ingredient in the proof is an estimate of the oscillation of the support function as it appears in \cite[Prop.~4, Lemma~5]{Andrews:/1999}, also cf. \cite[Prop.~7.3]{Schnuerer:02/2006}. These results hold for the case $n=2.$ In higher dimensions there is a generalization of these results, which provides closeness to a sphere in terms of the difference of the principal radii of a convex hypersurface, cf. \cite[Thm.~1.4]{Leichtweis:08/1999}.
 
In a recent paper Hung and Wang came up with a counterexample to such an asymptotical roundness for hypersurfaces expanding by the inverse mean curvature flow in the hyperbolic space, cf. \cite[Thm.~1]{HungWang:12/2014}. This shows the impossibility of proving results like ours in $\H^{n+1}.$

Before we give an overview over previous results on expanding flows in Euclidean space, let us state the main result of this paper. We require the following assumptions on $F.$

\Ass\label{F}
Let $\G\subset\R^{n}$ be an open, convex and symmetric cone containing the positive cone
\eq \G_{+}=\{(\k_{i})\in\R^{n}\cn \k_{i}>0,\ 1\leq i\leq n\}.\eeq
Let $F$ be a positive, monotone, symmetric and concave curvature function, normalized to $F(1,\dots,1)=n,$ such that
\begin{enumerate}[(i)]
\item{in case $0<p\leq 1$ we have $F\in C^{\infty}(\G)$ and $F_{|\del\G}=0,$}
\item{in case $p>1$ we additionally have $\G=\G_{+}.$}
\end{enumerate}
\eAss

Recall, that a hypersurface $M$ is called \textit{$F$-admissable}, if $F((\k_{i})(x))$ is well-defined for all $x\in M,$ where $\k_{i}(x)$ are the principal curvatures of $M$ at $x$ with respect to the inward unit normal.

The main result of this paper is the following one.

\Thm\label{Main}
Let $n\geq 2,$ $0<p<\infty$ and let $F$ satisfy Assumption \ref{F}. Let 
\eq x_{0}\cn M\hra M_{0}\subset\R^{n+1}\eeq
 be the smooth embedding of a closed, orientable, connected and $F$-admissable hypersurface, which can be written as a graph over a sphere $\S^{n},$
 \eq M_{0}=\{(u(0,x),x)\cn x\in \S^{n}\}.\eeq
Then
\begin{enumerate}[(i)]
\item\label{Maina}{there exists a unique smooth solution on a maximal time interval
\eq x\cn [0,T^{*})\times M\hra\R^{n+1},\eeq
which satisfies the flow equation
\SAl\label{Floweq} \dot{x}&=\fr{1}{F^{p}}\nu\\
		x(0,\xi)&=x_{0}(\xi),\end{split}\end{align}
where $\nu=\nu(t,\xi)$ is the outward unit normal to $M_{t}=x(t,M)$ at $x(t,\xi)$ and $F$ is evaluated at the principal curvatures of $M_{t}$ at $x(t,\xi).$}
\item\label{Mainb}{There exists a point $Q\in\R^{n+1}$ and a sphere $S^{*}=S_{R^{*}}(Q)$ around $Q$ with radius $R^{*},$ such that the spherical solutions $S_{t}$ with radii $R_{t}$ of (\ref{Floweq}) with $M_{0}=S_{R^{*}}$ satisfy
\eq \label{Maind} \mrm{dist}(M_{t},S_{t})\leq cR_{t}^{-\fr p2}\quad\forall t\in[0,T^{*}),\eeq $c=c(p,M_{0},F).$ Here $\mrm{dist}$ denotes the Hausdorff distance of compact sets.}
\end{enumerate}
\eThm

Statement (\ref{Maina}) is just the existence of a solution on a maximal time interval. This result is not new, holds in even more general situations and a proof can be found in \cite[Thm.~2.5.19, Lemma~2.6.1]{Gerhardt:/2006}. We stated it for convenience. Also note that the statement in (\ref{Mainb}) indeed says that the flow becomes close to a flow of spheres. This is due to the fact that the radii of spheres which satisfy (\ref{Floweq}) with a sphere as initial hypersurface do converge to infinity during the maximal time of existence, cf. \cite[Rem.~3.1]{Gerhardt:01/2014} and \cite[Thm.~0.1]{Gerhardt:/1990}.

Indeed, (\ref{Maind}) also allows to improve the rate of convergence of $\~{M}_{t}$ by choosing the \textit{optimal} geodesic sphere to rescale. We will not carry this out here, but refer to \cite[Sec.~7]{Schnuerer:02/2006} for a rough outline of the arguments involved.
 
Now we give a brief overview over the state of the art in classical expanding curvature flows. We leave aside the theory of contracting flows, weak solutions, flows with boundary conditions, other ambient spaces and evolving curves, due to the tremendous amount of literature, which is not of direct interest with respect to our results.

For smooth, expanding flows in Euclidean space usually, except for \cite{Schnuerer:/2006}, the asymptotic behavior of the flow hypersurfaces is described via the rescaling
\eq \~{M}_{t}=\T^{-1}M_{t},\eeq
where $\T(t)$ is the radius of the evolution of an \textit{arbitrary} geodesic sphere, such that $\~{M}_{t}$ is bounded below and above. In those works, the authors show that $\~{M}_{t}$ converges to a sphere smoothly, which is less than (\ref{Maind}) on the $C^{0}$-level.

  Convergence of the rescaled surfaces $\~{M}_{t}$ was proven, for example, in the papers by Gerhardt, \cite{Gerhardt:/1990}, and Urbas, \cite{Urbas:/1990}, in the case $p=1$ under similar assumptions on $F$ as we do impose them. Results like these in the case $0<p\leq 1$ and general $F$ were derived in \cite{Urbas:/1991} and for more general, but still concave functions of the inverse Gauss curvature or the principal radii, compare the works by Chow and Tsai, \cite{ChowTsai:/1998} and \cite{ChowTsai:/1997} respectively, as well as \cite{IvochkinaNehringTomi:/2000}. Results for $p>1$ have been accomplished for $n=p=2$ and $F=2K^{\fr 12},$ the classical inverse Gauss curvature flow in $\R^{3},$ by Schn{\"u}rer, \cite{Schnuerer:/2006}, and for $n=2,$ $1<p<2$ by Li, \cite{Li:06/2010}. Probably the most general existing paper on classical inverse curvature flows in $\R^{n+1}$ is \cite{Gerhardt:01/2014}. Besides those convergence results, Smoczyk has found an explicit representation for the solution of the inverse harmonic mean curvature flow, \cite{Smoczyk:/2005}.
  
To our knowledge, the only situation, in which statement (\ref{Mainb}) is proven, is the case $n=p=2$ and $F=2K^{\fr 12},$ where $K$ is the Gaussian curvature, cf. \cite{Schnuerer:/2006}. We are not aware of the existence of a convergence result of type (\ref{Maind}) in case of the other parameters.
 
%%%%%%%%%%%%%%%%%%%%%%%%%%%%%%%%%%%%%%%%%%%%%%%%%%%%%%%%%%%%%%%%%%%%%%%%%%%%%%%%
 
%%%%%%%%%%%%%%%%%%%%%%%%%%%%%%%%%%%%%%%%%%%%%%%%%%%%%%%%%%%%%%%%%%%%%%%%%%%%
\section{Notation and definitions}
In this article we consider closed, embedded and oriented hypersurfaces $M\hra\R^{n+1},$ which can be written as graphs over a sphere $\S^{n},$
\eq M=\{(u(x),x)\cn x\in\S^{n}\}.\eeq The coordinate representation $(u(x),x)$ is to be understood in polar coordinates, in which the Euclidean metric reads
\eq d\-{s}^{2}=dr^{2}+r^{2}\s_{ij}dx^{i}dx^{j}\equiv dr^{2}+\-{g}_{ij}dx^{i}dx^{j}.\eeq 
We prefer the coordinate based notation for tensors. 

Note that sometimes we use the slightly ambiguous notation to write $x$ for an element $x\in\S^{n},$ where it is to be understood as $x=(x^{i}),$ latin indices ranging between $1$ and $n,$ or to write $x$ as an element $x\in\R^{n+1},$ where it is then to be understood as $(x^{\a})=(x^{0},(x^{i})),$
greek indices ranging from $0$ to $n.$ Then $x^{0}$ denotes the radial component $r.$

We denote the standard induced metric of $\S^{n}\hra\R^{n+1}$ by $(\s_{ij}).$  The geometric quantities of $M$ are denoted via the following notation. The induced metric is denoted by $g=(g_{ij})$ with inverse $g^{-1}=(g^{ij})$ and the second fundamental form with respect to the inward normal is $A=(h_{ij}).$ Tensor indices of tensor fields on $M$ are always lowered or lifted via $g,$ unless stated otherwise, e.g.
\eq h^{i}_{j}=g^{ik}h_{kj}.\eeq
Covariant derivatives with respect to the induced metric will simply be denoted by indices, e.g. $u_{i}$ for a function $u\cn M\ra\R,$ or by a semicolon, if ambiguities are possible, e.g. $h_{ij;k}.$

The outward normal vector field to $M$ is given by
\eq (\nu^{\a})=v^{-1}(1,-\check{u}^{i}),\eeq
where $\check{u}^{i}=\-{g}^{ik}u_{k},$ $(\-{g}^{ik})=(\-{g}_{ik})^{-1}$ and 
\eq v^{2}=1+\-{g}^{ij}u_{i}u_{j}\equiv 1+|Du|^{2}.\eeq
For a tensor field $T=(t^{i_{1},\dots i_{k}}_{j_{1},\dots,j_{l}})$ on $M$ the pointwise norm $\|T\|$ is always defined with respect to the induced metric
\eq \|T\|^{2}=t^{i_{1},\dots,i_{k}}_{j_{1},\dots,j_{l}}t^{j_{1},\dots,j_{l}}_{i_{1},\dots,i_{k}}.\eeq 

A dot over a function or a tensor always denotes a total time derivative, e.g.
\eq \dot{u}=\fr{d}{dt}u,\eeq
whereas a prime denotes differentiation with respect to a direct argument. If for example $f=f(u),$ then 
\eq f'=\fr{d}{du}f\eeq
and \eq \dot{f}=f'\dot{u}.\eeq
Note, that this notation partially deviates from those in \cite{Gerhardt:/2006} and \cite{Gerhardt:01/2014}.

\TOCstop
\subsection*{Curvature functions}\
\TOCstart

The formulation of our assumptions on the curvature function $F$ has used cones in $\R^{n},$ i.e. $F$ depends smoothly on the principal curvatures,
\eq F=F(\k_{i}).\eeq However, as it is shown in \cite[Ch.~2.1]{Gerhardt:/2006} and the references therein, it is also possible to consider $F$ as a smooth function of the second fundamental form and the metric, 
\eq F=F(h_{ij}, g_{ij}),\eeq
or, as well, as a function defined on the mixed tensor $(h^{i}_{j}),$
\eq F=F(h^{i}_{j}).\eeq
Those formulations are basically equivalent. In the formulation of evolution equations we will use the second of those three. Note, that 
\eq F^{kl}=\fr{\del F}{\del h_{kl}}\eeq
defines a tensor field on $M$ of two contravariant indices.

\TOCstop
\subsection*{Evolution equations}\
\TOCstart

\Rem
The existence of a solution to (\ref{Floweq}) on a maximal time interval $[0,T^{*})$ is well-known. We refer to \cite[Thm.~2.5.19, Lemma~2.6.1]{Gerhardt:/2006}. Furthermore, the solution $x$ exists at least as long as the solution 
\eq u\cn [0,\-{T})\times \S^{n}\ra\R^{}\eeq of the scalar flow equation
\SAl \fr{\del}{\del t}u&=\fr{v}{F^{p}}\\
		u(0,\cdot)&=u_{0},\end{split}\end{align}
where $u_{0}$ is the graph representation of the initial hypersurface, also compare \cite[Thm.~2.5.17]{Gerhardt:/2006} and \cite[p.~98-99]{Gerhardt:/2006}. Note as well that under Assumption \ref{F} we have $\-{T}=T^{*},$ cf. \cite[Thm.~1.1, Thm.~1.2]{Gerhardt:01/2014}.
\eRem

For real numbers $r>0$ define
 \eq \Phi(r)=-r^{-p}.\eeq
The relevant evolution equations involved in the curvature flow are the following. %The metric $(g_{ij})$ satisfies, cf. \cite[Lemma 2.3.1]{Gerhardt:/2006},
%\eq \label{Evg} \dot{g}_{ij}=-2\Phi h_{ij},\eeq
%where $\Phi=\Phi(F).$
The second fundamental form %$(h_{ij})$ of the flow hypersurfaces evolves according to
%\SAl	\label{EvA}\dot{h}_{ij}-\Phi'F^{kl}h_{ij;kl}&=\Phi'F^{kl}h_{rk}h^{r}_{l}h_{ij}-(\Phi'F+\Phi)h_{ik}h^{k}_{j}+\Phi^{kl,rs}h_{kl;i}h_{rs;j}\end{split}\end{align}
%and
 in mixed form satisfies
\SAl	\dot{h}^{i}_{j}-\Phi'h^{i}_{j;kl}&=\Phi'F^{kl}h_{rk}h^{r}_{l}h^{i}_{j}-(\Phi'F-\Phi)h^{i}_{k}h^{k}_{j}+\Phi^{kl,rs}h_{kl;j}{h_{rs;}}^{i},\end{split}\end{align}
cf. \cite[Lemma~2.4.1]{Gerhardt:/2006}.
The curvature function $\Phi$ satisfies
\eq \dot{\Phi}-\Phi'F^{kl}\Phi_{kl}=\Phi'F^{kl}h_{rk}h^{r}_{l}\Phi,\eeq
cf. \cite[Lemma~2.3.4]{Gerhardt:/2006}.
In the sequel we will need two other derived evolution equations, namely for the mean curvature $H=h^{i}_{i},$
\eq \label{EvH}\dot{H}-\Phi'F^{kl}H_{kl}=\Phi'F^{kl}h_{rk}h^{r}_{l}H-(\Phi'F-\Phi)\|A\|^{2}+\Phi^{kl,rs}h_{kl;i}{h_{rs;}}^{i}\eeq
and for $\|A\|^{2}=h^{i}_{j}h^{j}_{i},$
\SAl \label{EvA2}	\fr{d}{dt}\|A\|^{2}-\Phi'F^{kl}(\|A\|)_{kl}&=2\Phi'F^{kl}h_{rk}h^{r}_{l}\|A\|^{2}-2(\Phi'F-\Phi)h^{i}_{k}h^{k}_{j}h^{j}_{i}\\
								&\hp{=}+2\Phi^{kl,rs}h_{kl;i}{h_{rs;}}^{j}h^{i}_{j}-2\Phi'F^{kl}h^{i}_{j;k}h^{j}_{i;l}.\end{split}\end{align}
												
\Rem
For better readability of the subsequent results, we stick to the convention that whenever we claim the existence of constants, $c,$ $\g$ etc., they are allowed and understood to depend on $n,$ $p,$ $M_{0}$ and $F$ as given data of the initial value problem, without mentioning this over and over again. 
\eRem		
%%%%%%%%%%%%%%%%%%%%%%%%%%%%%%%%%%%%%%%%%%%%%%%%%%%%%%%%%%%%%%%%%%%%%%%%%%%%%

%%%%%%%%%%%%%%%%%%%%%%%%%%%%%%%%%%%%%%%%%%%%%%%%%%%%%%%%%%%%%%%%%%%%%%%%%%%%%
\section{Pinching estimates}
In this section we successively improve the pinching estimates. First we need to revisit some results from \cite{Gerhardt:01/2014}.

Let $\T=\T(t,r)$ denote a geodesic sphere with initial radius $r,$ that exists as long as the solution $x$ of (\ref{Floweq}) and for which
\eq \~{u}=u\T^{-1}\eeq
is bounded below and above by positive constants, compare \cite[Lemma 3.3-3.5]{Gerhardt:01/2014} and \cite[(4.8)]{Gerhardt:/1990}. The following proposition holds.

\Prop\label{ScaledA}
Let $x$ be the solution of (\ref{Floweq}) under Assumption \ref{F}. Then the rescaled principal curvatures 
\eq \~{\k}_{i}=\k_{i}\Theta\eeq
satisfy
\eq 0<c^{-1}\leq\~{\k}_{i}\leq c\ \ \forall t\in[0,T^{*}),\eeq
if $p>1,$ and in case $p\leq 1$ the $\~{\k}_{i}$ stay within a compact subset of $\G.$ 
\eProp

\pf
For $p<1$ this is \cite[Lemma~4.11, Cor.~4.12]{Gerhardt:01/2014} and for $p=1$ we refer to \cite[(4.9)]{Gerhardt:/1990}. In case $p>1$ we refer to \cite[Lemmata~3.10, 4.7, 4.9]{Gerhardt:01/2014}.
\epf

Furthermore, we obtain decay estimates for the gradient.

\Prop\label{graddecay}
Let $x$ be the solution of (\ref{Floweq}) under Assumption \ref{F}, $0<p<\infty.$ Then there exist positive constants $c$ and $\g,$ such that the function
\eq \p=\log u\eeq
satisfies the gradient estimate
\eq |D\p|=\s^{ij}\p_{i}\p_{j}\leq c\T^{-\g}\ \ \forall t\in[0,T^{*}).\eeq
\eProp

\pf
In case $p>1$ this holds with $\g=\fr 12,$ cf. \cite[Lemma~3.7]{Gerhardt:01/2014}. In case $p=1$ this follows from \cite[Lemma~2.5]{Gerhardt:/1990}. In case $p<1$ this follows from the formula
\cite[(3.45)]{Gerhardt:01/2014}, where one should also note \cite[(3.41)]{Gerhardt:01/2014}, as well as the spherical growth
\eq\label{Theta} \T(t,r)=\(\fr{1-p}{n^{p}}t+r^{1-p}\)^{\fr{1}{1-p}},\eeq
cf. \cite[(3.11)]{Gerhardt:01/2014}.
\epf

From these observations we deduce that the pinching of the hypersurfaces actually improves.

\Prop\label{PinchB}
Let $x$ be the solution of (\ref{Floweq}) under Assumption \ref{F}, $0<p<\infty.$ Then there exist positive constants $c$ and $\g,$ such that the principal curvatures $\k_{i}$ of the flow hypersurfaces satisfy the pointwise estimate
\eq \label{PinchBa}	(\k_{i}-\k_{j})^{2}\leq cH^{2+\g}\ \ \forall t\in[0,T^{*}).\eeq
\eProp

\pf
For this proof norms of tensors are formed with respect to $\s_{ij}.$

The function $\log \~{u}$ is bounded in $C^{\infty}(\S^{n}),$ cf. \cite[Lemma~5.1, Thm.~5.2]{Gerhardt:01/2014}. Thus, via interpolation, compare \cite[Lemma~6.1]{Gerhardt:11/2011} and Proposition \ref{graddecay}, we obtain
\eq |D^{2}\p|^{2}\leq c|D\p||\p|_{C^{3}}\leq c\T^{-\g}\ \ \forall t\in[0,T^{*}).\eeq 
The rescaled second fundamental form $\~{h}^{i}_{j}=\T h^{i}_{j}$ satisfies
\eq \~{h}^{i}_{j}=v^{-1}\~{u}^{-1}(\d^{i}_{j}-(\s^{ik}-v^{-2}\p^{i}\p^{k})\p_{kj}),\eeq
cf. \cite[(5.2)]{Gerhardt:01/2014} and \cite[Lemma~2.7.6]{Gerhardt:/2006}.
Thus
\eq |\~{h}^{i}_{j}-\l\d^{i}_{j}|\leq c\T^{-\fr{\g}{2}},\eeq
where we used Proposition \ref{graddecay} and the fact that $\~u$ converges to some constant $\l^{-1}$. 
Rescaling backwards yields 
\eq |\k_{i}-\l\T^{-1}|\leq c\T^{-\(1+\fr{\g}{2}\)}\quad\forall 1\leq i\leq n,\eeq 
hence we obtain the result in view of
\eq 0<c^{-1}\leq \T H\leq c.\eeq
\epf

The final pinching improvement will allow us to derive asymptotical roundness of the flow hypersurfaces.

\Prop\label{PinchC}
Let $x$ be the solution of (\ref{Floweq}) under Assumption \ref{F}, $0<p<\infty.$ Then for all $\d<4+2p$ there exists a constant $c=c_{\d}>0,$ such that the principal curvatures $\k_{i}$ satisfy the pointwise estimate
\eq \label{PinchCa}	(\k_{i}-\k_{j})^{2}\leq c H^{\d}\quad \forall t\in[0,T^{*}).\eeq
\eProp

\pf
Using Proposition \ref{PinchB} it suffices to show this for $\d>2.$ From (\ref{EvH}) and (\ref{EvA2}) we obtain that
\eq w=\|A\|^{2}-\fr{1}{n}H^{2}\eeq
 satisfies
 \SAl	\dot{w}-\Phi'F^{kl}w_{kl}&=2\Phi'F^{kl}h_{rk}h^{r}_{l}w-2(\Phi'F-\Phi)\(h^{i}_{k}h^{k}_{j}h^{j}_{i}-\fr{1}{n}\|A\|^{2}H\)\\
 						&\hp{=}+2\Phi^{kl,rs}h_{kl;i}{h_{rs;}}^{j}\(h^{i}_{j}-\fr{1}{n}H\d^{i}_{j}\)\\
						&\hp{=}-2\Phi'F^{kl}\(h^{i}_{j;k}h^{j}_{i;l}-\fr 1n H_{k}H_{l}\)\end{split}\end{align}
and thus, for positive constants $c$ and $\d>2$ yet to be chosen, we obtain the evolution equation for 
\eq z=w-cH^{\d},\eeq
namely
\SAl\label{PinchC1}	\dot{z}-\Phi'F^{kl}z_{kl}&=2\Phi'F^{kl}h_{rk}h^{r}_{l}z-c(\d-2)\Phi'F^{kl}h_{rk}h^{r}_{l}H^{\d}\\
						&\hp{=}-2(\Phi'F-\Phi)\(h^{i}_{k}h^{k}_{j}h^{j}_{i}-\fr 1n\|A\|^{2}H-\fr{c\d}{2}H^{\d-1}\|A\|^{2}\)\\
						&\hp{=}+2\Phi^{kl,rs}h_{kl;i}{h_{rs;}}^{j}\(h^{i}_{j}-\fr 1n H\d^{i}_{j}-\fr{c\d}{2}H^{\d-1}\d^{i}_{j}\)\\
						&\hp{=}-2\Phi'F^{kl}\(h^{i}_{j;k}h^{j}_{i;l}-\fr 1n H_{k}H_{l}\)\\
						&\hp{=}+c\d(\d-1)H^{\d-2}\Phi'F^{kl}H_{k}H_{l}.\end{split}\end{align}
From Proposition \ref{PinchB} we find $0<t_{0}<T^{*},$ such that 
\eq \label{PinchC2}  \max\(\fr{w}{H^{2}},\k_{n}\)<\s:=\fr{1}{\d^{\a}n(n-1)}\ \ \forall t\in[t_{0},T^{*}),\eeq
where $\a_{\d}>> 1$ will be chosen appropriately later. Thus $t_{0}$ will only depend on $\d,$ $p$ and $M_{0}$ as well and thus we may define
\eq c=\fr{\s}{\inf\limits_{M_{t_{0}}}H^{\d-2}}.\eeq
Note that due to \cite[Lemma~2.2]{AndrewsMcCoy:03/2012} the $M_{t}$ are strictly convex for $t\geq t_{0}.$
On $M_{t_{0}}$ we have
\eq z=w-cH^{\d}=H^{2}\(\fr{w}{H^{2}}-cH^{\d-2}\)<0.\eeq
We wish to show that this remains valid up to $T^{*}.$ Thus supposeÊ $t_{1}>t_{0}$ to be the first time, such that there exists $\xi_{1}\in M_{t_{1}}$ with the property
\eq z(t_{1},\xi_{1})=0.\eeq
Define
\eq \e:=cH^{\d-2}(t_{1},\xi_{1}),\eeq
then there holds
\eq 0<\e<\s,\eeq
due to (\ref{PinchC2}).
From (\ref{PinchC1}) we obtain at $(t_{1},\xi_{1})$ that
\SAl \label{PinchC3}	0&\leq -c(\d-2)\Phi'F^{kl}h_{rk}h^{r}_{l}H^{\d}\\
		&\hp{=}-2(\Phi'F-\Phi)\(h^{i}_{k}h^{k}_{j}h^{j}_{i}-\fr1n \|A\|^{2}H-\fr{c\d}{2}H^{\d-1}\|A\|^{2}\)\\
		&\hp{=}+2\Phi^{kl,rs}h_{kl;i}{h_{rs;}}^{j}\(h^{i}_{j}-\fr 1nH\d^{i}_{j}-\fr{c\d}{2}H^{\d-1}\d^{i}_{j}\)\\
		&\hp{=}-2\Phi'F^{kl}\(h^{i}_{j;k}h^{j}_{i;l}-\fr 1n H_{k}H_{l}\)+c\d(\d-1)H^{\d-2}\Phi'F^{kl}H_{k}H_{l}.\end{split}\end{align}
From \cite[Lemma~2.3]{AndrewsMcCoy:03/2012} we have at $(t_{1},\xi_{1}),$ note that $w=\e H^{2},$
\SAl h^{i}_{k}h^{k}_{j}h^{j}_{i}-\(\fr 1n+\e\)H\|A\|^{2}&\geq \e\(\fr{1}{n}+\e\)(1-\sqrt{n(n-1)\e})H^{3}\\
								&=\e\(1-\sqrt{n(n-1)\e}\)\|A\|^{2}H\\
								&>0\end{split}\end{align}
and we have
\SAl	g^{kl}\(h^{i}_{j;k}h^{j}_{i;l}-\fr 1n H_{k}H_{l}\)&=\left\|D\(A-\fr 1n Hg\)\right\|^{2}\\
									&\geq \fr{2(n-1)}{3n}\|DA\|^{2}\\
									&\geq \fr{2(n-1)}{n(n+2)}\|DH\|^{2},\end{split}\end{align}
cf. \cite[Lemma~2.1]{AndrewsMcCoy:03/2012}. In view of the concavity of $F$ we have $F\leq H,$ \cite[Lemma~2.2.20]{Gerhardt:/2006}, and thus from (\ref{PinchC3}) we obtain
\SAl	0&\leq -2(\Phi'F-\Phi)\(\e\(2-\sqrt{n(n-1)\e}\)\|A\|^{2}H-\fr{\e\d}{2}\|A\|^{2}H\)\\
		&\hp{=}-\e(\d-2)\fr{p}{p+1}(\Phi'F-\Phi)\|A\|^{2}H\\
		&\hp{=}-\e(\d-2)\fr{p}{p+1}(\Phi'F-\Phi)(F^{kl}-g^{kl})h_{rk}h^{r}_{l}H\\
		&\hp{=}+2H\Phi^{kl,rs}h_{kl;i}{h_{rs;}}^{j}\(H^{-1}h^{i}_{j}-\fr 1n \d^{i}_{j}-\fr{\e\d}{2}\d^{i}_{j}\)-\fr{2(n-1)}{n(n+2)}\Phi'\|DH\|^{2}\\
		&\hp{=}-\fr{2(n-1)}{3n}\Phi'\|DA\|^{2}-2\Phi'(F^{kl}-g^{kl})\(h^{i}_{j;k}h^{j}_{i;l}-\fr 1n H_{k}H_{l}\)\\
		&\hp{=}+\e\d(\d-1)\Phi'\|DH\|^{2}+\e\d(\d-1)\Phi'(F^{kl}-g^{kl})H_{k}H_{l}. \end{split}\end{align} 
In view of (\ref{PinchC2}) and due to to the fact that 
\eq \|F^{kl}-g^{kl}\|\ra 0,\eeq
we may without loss of generality enlarge $t_{0}$ and $\a,$ such that the terms involving curvature derivatives are absorbed by the terms
\eq -\fr{2(n-1)}{n(n+2)}\Phi'\|DH\|^{2}\ \mrm{and}\ -\fr{2(n-1)}{3n}\Phi'\|DA\|^{2}.\eeq
In this case at $(t_{1},\xi_{1})$
 we obtain \SAl\label{PinchC4}		0&\leq -2(\Phi'F-\Phi)\e\|A\|^{2}H\(\fr{p+2}{p+1}-\d^{-\fr{\a}{2}}-\fr{\d}{2(p+1)}\)\\
			&\hp{=}-\e(\d-2)\fr{p}{p+1}(\Phi'F-\Phi)(F^{kl}-g^{kl})h_{rk}h^{r}_{l}H\\
			&<0,\end{split}\end{align}
if we choose
\eq \d<4+2p,\eeq
$\a=\a(p,\d)$ large enough and again $t_{0}$ larger to ensure that the term involving $\|F^{kl}-g^{kl}\|$ can be absorbed by the strictly negative, remaining part of the first line
in (\ref{PinchC4}).
This contradiction shows that $z$ will remain negative up to the time $T^{*},$ which yields the result.
\epf
%%%%%%%%%%%%%%%%%%%%%%%%%%%%%%%%%%%%%%%%%%%%%%%%%%%%%%%%%%%%%%%%%%%%%%%%%%%%%

%%%%%%%%%%%%%%%%%%%%%%%%%%%%%%%%%%%%%%%%%%%%%%%%%%%%%%%%%%%%%%%%%%%%%%%%%%%%%%
\section{Oscillation decay}\label{OscDecay}
Let $\rho_{+}(t)$ denote the \textit{circumradius} of $M_{t},$ i.e. the radius of the smallest ball in $\R^{n+1}$ enclosing $M_{t}$ and analogously $\rho_{-}(t)$ the \textit{inradius} of $M_{t},$ the radius of the largest ball in $\R^{n+1}$ enclosed by $M_{t}.$ In this section we prove (\ref{Mainb}) of Theorem \ref{Main}, namely that $\rho_{+}-\rho_{-}$ converges to zero and that we find an expanding family of geodesic spheres $S_{t}$ with radii $R_{t}$ the flow hypersurfaces $M_{t}$ fit themselves to,
\eq \dist(M_{t},S_{t})<cR_{t}^{-\fr p2}.\eeq
 Note that the $R_{t}$ represent a suitable choice of the $\T(t).$

Let $\hat{M}_{t}$ denote the convex body enclosed by $M_{t},$ which is well defined for large $t.$ For an interior point $y\in\mrm{int}(M_{t})$ let $u_{y}$ denote the graph representation of $M_{t}$ over the standard sphere centered in $y,$
\eq M_{t}=\{(u_{y}(x),x)\cn x\in \S^{n}(y)\}.\eeq

\Prop\label{MassCenter}
Let $x$ be the solution of (\ref{Floweq}) under Assumption \ref{F}, $0<p<\infty.$ Then there holds
\eq \rho_{+}(t)-\rho_{-}(t)\leq \osc u_{y_{t}}\leq c\T^{-\fr p2}(t)\ \ \forall t\in[0,T^{*}),\eeq
where $y_{t}$ is a suitable oscillation minimizing center of $\hat M_{t}.$
\eProp

\pf
First note that $\T$ is still an arbitrary rescaling factor as in Proposition \ref{ScaledA}.

From (\ref{PinchCa}) we deduce that principal curvatures satisfy
\eq ({\k}_{i}-{\k}_{j})^{2}\leq cH^{\d}\leq c\T^{-\d}\eeq
for $\d=4+p.$ Since for large $t$ the $M_{t}$ are strictly convex we may consider the difference of the largest and smallest principal radius of curvature,
\eq \left|\fr{1}{\k_{1}}-\fr{1}{\k_{n}}\right|=\left|\fr{\k_{n}-\k_{1}}{\k_{1}\k_{n}}\right|\leq c\T^{2-\fr{\d}{2}}.\eeq
Applying \cite[Thm.~1.4]{Leichtweis:08/1999} we obtain
\eq \dist(M_{t},S(y_{t}))\leq c\T^{2-\fr{\d}{2}},\eeq
from which the claim follows. Note that if a $y_{t}$ happens not to minimize the oscillation, we can adjust it to do so. A simple geometric argument also shows that then the $y_{t}$ must lie in the convex body $\hat M_{t}.$ 
\epf

In order to find an optimally fitted spherical flow, we need the centers $y_{t}$ from Proposition \ref{MassCenter} to converge.

\Lem\label{Q}
The centers $y_{t}$ from Proposition \ref{MassCenter} converge in $\R^{n+1}.$ In particular we have
\eq |y_{t}-Q|\leq c\T^{-\fr p2}(t)\eeq
for the limit $Q\in\R^{n+1}.$
\eLem

\pf
%This is a standard argument, using that the oscillation of the support function with respect to a fixed point is decreasing, cf. \cite[Thm.~3.1]{McCoy:/2003}, originally proved by Chow and Gulliver, \cite{ChowGulliver:/1996}. The argument is performed in detail in \cite[Lemma~7.5]{Schnuerer:02/2006}.
For a point $q\in\R^{n+1}$ let 
\eq\-u_{q}=\inpr{x-q}{\nu}\eeq
denote the support function with respect to $q.$ For $t$ close to $T^{*}$ the $M_{t}$ are strictly convex and hence we may apply a gradient estimate for convex hypersurfaces, \cite[Lemma~2.7.10]{Gerhardt:/2006}, to conclude
\eq v\leq e^{\-\k\osc u_{y_{t}}},\eeq
where $\-\k$ is an upper for the principal curvatures of the slices $\{x^{0}=\mrm{const}\},$
which in our case can be estimated:
\eq \-\k\leq \max_{M_{t}}\fr{1}{u_{y_{t}}}.\eeq
Thus
\eq v-1\leq e^{c\-\k\T^{-\fr{p}{2}}}-1\leq c\-\k \T^{-\fr{p}{2}}.\eeq
Since 
\eq \fr{u_{y_{t}}}{v}=\-u_{y_{t}},\eeq
we obtain
\eq |\-u_{y_{t}}-u_{y_{t}}|\leq u_{y_{t}}\fr{v-1}{v}\leq c\T^{-\fr{p}{2}}\eeq
and hence
\eq\max \-u_{y_{t}}-\min\-u_{y_{t}}\leq\max u_{y_{t}}-\min u_{y_{t}}+c\T^{-\fr{p}{2}}\leq c\T^{-\fr{p}{2}}(t).\eeq
Using that the oscillation of the support function with respect to a fixed point is decreasing, cf. \cite[Thm.~3.1]{McCoy:/2003}, originally proved by Chow and Gulliver, \cite{ChowGulliver:/1996}, we obtain
\eq \osc u_{y_{t_{2}}}\leq \osc\-u_{y_{t_{2}}}\leq \osc \-u_{y_{t_{1}}}\leq c\T^{-\fr p2}(t_{1})\quad\forall t_{2}>t_{1}\eeq
and  since $y_{t_{2}}$ is oscillation minimizing we must have
\eq |y_{t_{1}}-y_{t_{2}}|\leq c\T^{-\fr p2}(t_{1})\quad \forall t_{2}>t_{1}.\eeq
Letting $t_{1}\ra T^{*}$ we obtain the limit $Q$ and then for fixed $t_{1}$ letting $t_{2}\ra T^{*},$ we obtain the desired estimate.
\epf

We finish this paper by showing that the flow actually becomes close to a flow of spheres.

\Thm\label{Final}
Let $x$ be the solution of (\ref{Floweq}) under Assumption \ref{F}. Then there exists a geodesic sphere $S_{R^{*}}(Q),$ where $Q$ is the limit from Lemma \ref{Q}, such that the spherical leaves $S_{t}$ of the initial value problem
\SAl	\dot{y}&=\fr{1}{F^{p}}\nu\\
		y(0,M)&=S_{R^{*}}(Q)\end{split}\end{align}
and the flow hypersurfaces $M_{t}$ of the flow $x$ satisfy
\eq \dist (M_{t},S_{t})<cR_{t}^{-\fr p2},\eeq
where $R_{t}$ is the radius of $S_{t}.$		
\eThm

\pf
In case $p>1,$ $R^{*}$ is determined by the requirement, that the spherical flow exists as long as the flow $x.$ Then for all large $t$ we must have
\eq\label{Finala} S_{t}\cap M_{t}\neq \emptyset,\eeq
compare the arguments at the end of the proof of \cite[Lemma~5.1]{Schnuerer:/2006}. The result follows from the propositions \ref{MassCenter} and \ref{Q}.

In case $p\leq 1$ we need a different argument to show that there exists an expanding flow of spheres with the property (\ref{Finala}), since we have $T^{*}=\infty$ and $R^{*}$ is not determined clearly. Recall (\ref{Theta}), that for initial radius $r$ the radius $R$ of a sphere evolves according to
\eq R(t,r)=\(\fr{1-p}{n^{p}}t+r^{1-p}\)^{\fr{1}{1-p}},\eeq
if $p< 1$ and 
\eq R(t,r)=re^{\fr{t}{n}},\eeq
if $p=1.$ We see that in any case
\eq R_{t}=R(t,\cdot)\eeq
is an increasing diffeomorphism from $(0,\infty)$ onto its image $\(\(\fr{1-p}{n^{p}}t\)^{\fr{1}{1-p}},\infty\)$ in case $p<1$ and onto $(0,\infty)$ in case $p=1.$
Now let $u=u_{Q}$ and define sequences
\eq \-{R}^{k}=R_{k}^{-1}(\sup u(k,\cdot)),\eeq
\eq \bar{R}_{k}=R^{-1}_{k}(\inf u(k,\cdot))\eeq
and 
\eq R^{k}=\fr 12(\-{R}^{k}+\bar{R}_{k}).\eeq
By the maximum principle $\-R^{k}$ is non-increasing.
There holds
\eq S_{k}^{k}\cap M_{k}\neq\emptyset,\eeq
where $S^{k}_{k}$ denotes the spherical leave at time $k,$ which has started with initial radius $R^{k}.$
Since in case $p<1$
\eq \fr{d}{dr}R_{k}(r)=\(\fr{1-p}{n^{p}}k+r^{1-p}\)^{\fr{p}{1-p}}r^{-p}\eeq
and in case $p=1$
\eq \fr{d}{dr}R_{k}(r)=e^{\fr kn},\eeq
$\fr{d}{dr} R_{k}$ is uniformly bounded from below and since
 \eq \osc u\ra0,\eeq
 we obtain that $\-{R}^{k},$ $\-{R}_{k}$ and $R^{k}$ all converge to the same limit $R^{*}.$ We claim, that the initial sphere $S_{R^{*}}$ around $Q$ leads to a spherical flow satisfying (\ref{Finala}). 
  
 Otherwise there existed a time $k_{0},$ such that without loss of generality
 \eq R(k_{0},R^{*})<\inf u(k_{0},\cdot).\eeq
 By continuity of ODE orbits with respect to initial values on compact intervals, there is $\~{R}>R^{*},$ such that
 \eq R(k,R^{*})<R(k,\~{R})<\inf u(k,\cdot)\ \ \forall k\geq k_{0},\eeq
 where we also used the maximum principle.
 Applying $R_{k}^{-1}$ we find
 \eq R^{*}<\~{R}<\-{R}_{k}\ \ \forall k\geq k_{0},\eeq
 since $R^{-1}_{k}$ is increasing with respect to $r.$ This is a contradiction to $\-{R}_{k}\ra R^{*}.$ So the spherical leaves of the flow with initial value $S_{R^{*}}(Q)$ intersect the $M_{t}$ for all large times and due to the oscillation estimates we obtain the desired result.
\epf
%%%%%%%%%%%%%%%%%%%%%%%%%%%%%%%%%%%%%%%%%%%%%%%%%%%%%%%%%%%%%%%%%%%%%%%%%%%%%%%%

%%%%%%%%%%%%%%%%%%%%%%%%%%%%%%%%%%%%%%%%%%%%%%%%%%%%%%%%%%%%%%%%%%%%%%%%%%%%%%%%%
\section{Concluding remarks}
The pinching estimates, Proposition~\ref{PinchC}, turned out to improve, whenever $p$ becomes larger. This fact is somehow surprising, since the evolution equation of the gradient function
\eq v^{2}=1+|Du|^{2}\eeq
does not allow to apply the classical maximum principle, also compare the proof of \cite[Lemma~3.6]{Gerhardt:01/2014}, so one could expect that it should become harder to control oscillations. As we have seen, however, the equation for the traceless second fundamental form serves as a way out.

Also note that in a further work we applied a method similar to the one in section \ref{OscDecay} to prove that there can not be an estimate of the form
\eq\mrm{dist}(M,S_{R})\leq c\|\mathring{A}\|^{\a},\quad \a>1,\eeq
in the class of uniformly convex hypersurfaces with a universal constant. The idea is that otherwise we could use a similar proof as in section \ref{OscDecay} to prove asymptotical roundness of the inverse mean curvature flow in the hyperbolic space, which has shown not to be true in \cite{HungWang:12/2014}. See \cite{RothScheuer:05/2015} for a preprint version and a detailed description of this result.

\bibliographystyle{/Users/J_Mac/Documents/Uni/TexTemplates/hamsplain}
\bibliography{/Users/J_Mac/Documents/Uni/TexTemplates/Bibliography}

\providecommand{\bysame}{\leavevmode\hbox to3em{\hrulefill}\thinspace}
\providecommand{\href}[2]{#2}
\begin{thebibliography}{10}

\bibitem{Andrews:/1999}
Ben Andrews, \emph{Gauss curvature flow: the fate of the rolling stones},
  Invent. Math. \textbf{138} (1999), no.~1, 151--161.

\bibitem{AndrewsMcCoy:03/2012}
Ben Andrews and James McCoy, \emph{Convex hypersurfaces with pinched principal
  curvatures and flow of convex hypersurfaces by high powers of curvature},
  Trans. Am. Math. Soc. \textbf{364} (2012), no.~7, 3427--3447.

\bibitem{ChowGulliver:/1996}
Bennett Chow and Robert Gulliver, \emph{Aleksandrov reflection and nonlinear
  evolution equations, i: The n-sphere and n-ball}, Calc. Var. Partial Differ.
  Equ. \textbf{4} (1996), no.~3, 249--264.

\bibitem{ChowTsai:/1997}
Bennett Chow and Dong-Ho Tsai, \emph{Expansion of convex hypersurfaces by
  nonhomogeneous functions of curvature}, Asian J. Math. \textbf{1} (1997),
  no.~4, 769--784.

\bibitem{ChowTsai:/1998}
\bysame, \emph{Nonhomogeneous {G}auss curvature flows}, Indiana Univ. Math. J.
  \textbf{47} (1998), no.~3, 965--994.

\bibitem{Gerhardt:/1990}
Claus Gerhardt, \emph{Flow of nonconvex hypersurfaces into spheres}, J. Differ.
  Geom. \textbf{32} (1990), no.~1, 299--314.

\bibitem{Gerhardt:/2006}
\bysame, \emph{Curvature problems}, Series in Geometry and Topology, vol.~39,
  International Press of Boston Inc., 2006.

\bibitem{Gerhardt:11/2011}
\bysame, \emph{Inverse curvature flows in hyperbolic space}, J. Differ. Geom.
  \textbf{89} (2011), no.~3, 487--527.

\bibitem{Gerhardt:01/2014}
\bysame, \emph{Non-scale-invariant inverse curvature flows in {E}uclidean
  space}, Calc. Var. Partial Differ. Equ. \textbf{49} (2014), no.~1-2,
  471--489.

\bibitem{HungWang:12/2014}
Pei-Ken Hung and Mu~Tao Wang, \emph{Inverse mean curvature flows in the
  hyperbolic 3-space revisited}, Calc. Var. Partial Differ. Equ. (2014),
  {DOI:10.1007/s00526-014-0780-3}.

\bibitem{IvochkinaNehringTomi:/2000}
Nina Ivochkina, Thomas Nehring, and Friedrich Tomi, \emph{Evolution of
  starshaped hypersurfaces by nonhomogeneous curvature functions}, Algebra i
  Anal. \textbf{12} (2000), no.~1, 185--203.

\bibitem{Leichtweis:08/1999}
Kurt Leichtwei\ss, \emph{Nearly umbilical ovaloids in the n-space are close to
  spheres}, Result. Math. \textbf{36} (1999), no.~1-2, 102--109.

\bibitem{Li:06/2010}
Qi-Rui Li, \emph{Surfaces expanding by the power of the {G}auss curvature
  flow}, Proc. Amer. Math. Soc. \textbf{138} (2010), no.~11, 4089--4102.

\bibitem{McCoy:/2003}
James McCoy, \emph{The surface area preserving mean curvature flow}, Asian J.
  Math. \textbf{7} (2003), no.~1, 7--30.

\bibitem{RothScheuer:05/2015}
Julien Roth and Julian Scheuer, \emph{Explicit rigidity of almost-umbilical
  hypersurfaces}, preprint available at
  \href{http://arxiv.org/abs/1504.05749}{arxiv:1504.05749}, 2015.

\bibitem{Schnuerer:/2006}
Oliver~C. Schn{\"u}rer, \emph{Surfaces expanding by the inverse {G}au{\ss}
  curvature flow}, J. Reine Angew. Math. \textbf{600} (2006), 117--134.

\bibitem{Schnuerer:02/2006}
\bysame, \emph{Surfaces expanding by the inverse {G}auss curvature flow}, 2006,
  \href{http://arxiv.org/abs/math/0412297}{arxiv:0412297v2}.

\bibitem{Smoczyk:/2005}
Knut Smoczyk, \emph{A representation formula for the inverse harmonic mean
  curvature flow}, Elem. Math. \textbf{60} (2005), no.~2, 57--65.

\bibitem{Urbas:/1990}
John Urbas, \emph{On the expansion of starshaped hypersurfaces by symmetric
  functions of their principal curvatures}, Math. Z. \textbf{205} (1990),
  no.~1, 355--372.

\bibitem{Urbas:/1991}
\bysame, \emph{An expansion of convex hypersurfaces}, J. Differ. Geom.
  \textbf{33} (1991), no.~1, 91--125.

\end{thebibliography}

\end{document}